\documentclass[11pt, a4paper, reqno]{amsart}

\usepackage[pdftitle={On two notions of curvature on singular surfaces},
        pdfauthor={Maxime Marot}
        ]{hyperref}
\usepackage{amsmath, amsthm, amssymb, amsrefs, mathrsfs}
\usepackage{esint}
\usepackage{enumitem}
\usepackage[reftex]{theoremref}
\usepackage[skip=1mm, indent=15pt]{parskip}
\usepackage[left=2.5cm, right=2.5cm]{geometry}

\theoremstyle{definition}
\newtheorem{definition}{Definition}

\theoremstyle{plain}
\newtheorem{proposition}{Proposition}
\newtheorem{theorem}{Theorem}
\newtheorem{lemma}{Lemma}
\newtheorem{corollary}{Corollary}

\newtheorem*{corollary*}{Corollary}

\newcommand{\dd}{\ensuremath{\,\mathrm{d}}}
\newcommand{\vol}{\ensuremath{\mathrm{vol}}}
\newcommand{\N}{\ensuremath{\mathbb{N}}}
\newcommand{\Z}{\ensuremath{\mathbb{Z}}}
\newcommand{\R}{\ensuremath{\mathbb{R}}}
\newcommand{\C}{\ensuremath{\mathbb{C}}}
\newcommand{\loc}{\ensuremath{\mathrm{loc}}}

\DeclareMathOperator{\diam}{diam}

\DeclareMathOperator{\supp}{supp}
\DeclareMathOperator{\modang}{\tilde\measuredangle}
\DeclareMathOperator{\CBB}{CBB}
\DeclareMathOperator{\CAT}{CAT}

\DeclareMathOperator{\BICB}{BICB}
\DeclareMathOperator{\BICA}{BICA}
\DeclareMathOperator{\RCD}{RCD}

\title{On two notions of curvature on singular surfaces}

\author{Maxime Marot}

\address{Mathematisches Seminar, Heinrich-Hecht-Platz 6, 24118 Kiel, Germany}
\email{marot@math.uni-kiel.de}

\begin{document}

\begin{abstract}
In this paper, we investigate the equivalence of two distinct notions of curvature bounds on singular surfaces.
The first notion involves inequalities of the form $\omega\geq\kappa\mu$ (resp. $\omega\leq\kappa\mu$) where
$\omega$ is the curvature measure and $\mu$ the Hausdorff measure.
The second notion is the classical Alexandrov curvature bound CBB (resp. CAT).
We demonstrate that these two definitions are, in fact, equivalent.
Specifically, we fill an important gap in the theory by showing that the inequalities imply the corresponding Alexandrov CBB (resp. CAT) bound.
One striking application of our result is that, in combination with a result of Petrunin,
the lower bound $\omega\geq\kappa\mu$ implies $\RCD(\kappa, 2)$.
\end{abstract}

\maketitle

\section{Introduction}

In the second half of the twentieth century, the Leningrad school and A. D. Alexandrov developed synthetic notions of curvature bounds on metric spaces.
One of these is the well-known curvature bounded below (CBB) formalized in the seminal paper \cite{buragoADAlexandrovSpaces1992}
by Yu. Burago, M. Gromov, and G. Perel'man.
At the same time, Alexandrov and Zalgaller introduced the theory of surfaces of bounded integral curvature (BIC)
in the book \cite{aleksandrovIntrinsicGeometrySurfaces2002}.
The theory of CBB/CAT Alexandrov spaces is well-developed and continues to be a highly active area of research.
The objective of this paper is to demonstrate that, under a suitable definition of curvature bounds for BIC surfaces, the two theories are equivalent.

Let first define the notion of curvature bound:

\begin{definition}
Let $\kappa\in\mathbb{R}$.
We say that a BIC surface $S$ has \textit{curvature bounded below by $\kappa$}, shortly $S$ is $\BICB(\kappa)$,
if the inequality $\omega\geq\kappa\mu$ holds in the sense of set functions, i.e. for every Borel set $E\subset S$,
the inequality $\omega(E)\geq\kappa\mu(E)$ holds.
Analogously, we say that $S$ has \textit{curvature bounded above by $\kappa$}, shortly $S$ is $\BICA(\kappa)$,
if the inequality $\omega\leq\kappa\mu$ holds in the sense of set functions.
\end{definition}

One can easily see that it generalizes lower and upper bounds on Gaussian curvature on smooth Riemannian manifolds.

Our first main result is the following:

\begin{theorem}
\th\label{thm:BICBImpliesCBB}
Assume $S$ is a complete $\BICB(\kappa)$ surface for some $\kappa\in\R$
and $S$ has no cusp (there are no point $p\in S$ with $\omega(\{p\})\geq2\pi$).
Then $S$ is $\CBB(\kappa)$.
\end{theorem}

One can actually drop the assumption of absence of cusp when $\kappa\geq 0$.
The converse of this theorem is well-known in the literature:
if $S$ is a $\CBB(\kappa)$ surface of Hausdorff dimension $2$ then it is a BIC surface without any cusp.
see \cite{thomasrichardCanonicalSmoothing2018} for a geometric proof and
\cite{ambrosioRegularityAlexandrovSurfaces2016} for an analytic proof.

From Petrunin's work \cite{petruninAlexandrovMeetsLottVillaniSturm2011}, $\CBB(\kappa)$ spaces together with the
2-dimensional Hausdorff measure satisfy the $\RCD(\kappa, 2)$ condition.
Hence, an immediate and very important corollary is:

\begin{corollary}
\th\label{th:BICBImpliesRCD}
Assume $S$ is a complete BIC surface whose curvature measure satisfies $\omega\geq\kappa\mu$ for some $\kappa\in\R$.
Moreover, when $\kappa<0$, $S$ is supposed to be whithout cusp.
Then S together with its 2-dimensional Hausdorff measure is $\RCD(\kappa,2)$.
\end{corollary}

We also establish an analog of Theorem \ref{thm:BICBImpliesCBB} for curvature bounded above by $\kappa$,
i.e. when $\omega\leq\kappa\mu$.
The second main result is the following:

\begin{theorem}
\th\label{thm:BICAImpliesCAT}
Assume $S$ is a complete BIC surface whose curvature measure satisfies $\omega\leq\kappa\mu$ for some $\kappa<0$.
Then $S$ is locally $\CAT(\kappa)$.
\end{theorem}

The converse of this theorem has been shown recently in the paper \cite{chowdhuryCATKSurfaces2025}.
The conclusion of locally CAT is sharp.
For example, if one take a flat cylinder and a triangle on a circular section with all three sides with same length then the
angles are $\pi$.
Whereas the angles of a plane triangle with same side length has angles equal to $\pi/3$.
So the flat cylinder cannot be globally $\CAT(0)$ even though it has $\omega=0$.

Theorem \ref{thm:CmBBDiscCvxImpliesHinge} and Theorem \ref{thm:CmBADiscImpliesHinge}
have already been noted by Reshetnyak in \cite[p.140]{reshetnyakGeometryIV1993},
however, the results were only proved for $\kappa=0$ and in a rather sketchy manner.
Also we were not able to tell wether the theorems in Reshetnyak's survey are stated for upper angle
because the proofs Reshetnyak refers to are clearly for sector angles which are greater than upper angles.

The paper is structured as follow:
we first prove \thref{thm:BICBImpliesCBB} by showing that on every isothermal chart the non-smooth metric can approximate
by carefully crafted smooth metrics with curvature bounded below uniformly by $\kappa$.
For positive curvature the construction is a straightforward mollification but for negative curvature bound the proof is more involved.
One first need to mollify the curvature and then solve a perturbed Liouville-type equation prescribing the curvature for the conformal factor.
The technically involved part lies in showing that the curvature measures of these smooth metrics weakly converge to the original one.
Then the conclusion follows by Reshetnyak's convergence theorem and CBB globalization theorem.
In the second part, we discuss a synthetic approach to this problem and its limitation due to convexity of triangles and Gauss-Bonnet theorem.
In the case of curvature bounded above, this approach still leads to a proof of \thref{thm:BICAImpliesCAT}.

Note that completeness is only required for the globalization and one cannot expect to use \cite{petruninGlobalization2016}
-- which only requires geodesic distance but concludes on the closure -- as the completion of a BIC surface is unpredictable
(see \cite[Remark p.205]{aleksandrovIntrinsicGeometrySurfaces2002}).

\subsection*{Acknowledgement}

I would like to express my gratitude to Sebastian Boldt and my advisor Batu Güneysu for their precious help on tackling this problem
and pointing mistakes I made in my first tries.
This work would not have been possible without their guidance.
I also thank François Fillastre for answering many of my questions on BIC surfaces.

\section{Background}

\subsection{Angles in metric spaces}

We denote by $\mathbb{M}^2(\kappa)$, called \textit{model spaces}, the connected two-dimensional Riemannian manifold of constant sectional
curvature $\kappa$ and write $\varpi^\kappa$ its diameter.
Let $X$ be a metric space with distance function $d$.
For any $\kappa\in\R$ and any points $p,q,r\in X$, we define the \textit{model triangle} $\tilde\Delta^\kappa(pqr)$
to be any triangle $[\tilde p\tilde q\tilde r]$ in the model space $\mathbb{M}^2(\kappa)$ such that
\[
d_{\mathbb{M}^2(\kappa)}(\tilde p, \tilde q) = d(p,q),\quad
d_{\mathbb{M}^2(\kappa)}(\tilde q, \tilde r) = d(q,r),\quad
d_{\mathbb{M}^2(\kappa)}(\tilde r, \tilde p) = d(r,p).
\]
Note that when $\kappa\leq 0$, the model triangle always exists.
Otherwise, when $\kappa$ is positive, we require that the following condition holds,
\[
d(p,q) + d(q,r) + d(r,p) < 2 \varpi^\kappa.
\]
This construction allows us to define the \textit{model angle}, written $\modang^\kappa(p^q_r)$, by the angle at the vertex $\tilde p$
in the model triangle $[\tilde p\tilde q\tilde r]$.
Let $p,x,y\in X$ be a triple of points such that $p$ is different from $x$ and $y$.
A pair of minimizing geodesics $(\gamma, \gamma')$ such that $\gamma$ (resp. $\gamma'$) is from $p$ to $x$ (resp. $y$)
is called a \textit{hinge} and is denoted $[p^x_y]$.
Given a hinge $[p^x_y]$, we define its \textit{(upper) angle} to be
\[
\measuredangle[p^x_y] := \limsup_{\substack{q,r\to p \\ q,r\neq p}} \modang^\kappa(p^q_r).
\]
This definition doesn't depend on the $\kappa$ chosen.
We say that this angle is \textit{defined} if the $\limsup$ is actually a limit.
Two hinges $[p^x_y]$ and $[p^y_z]$ are said to be \textit{adjacent} if the union of the geodesics from $p$ to $x$ and from $p$ to $z$
is actually a geodesic from $x$ to $z$.

Another notion of angle is defined as follow.
Consider a hinge $[p^x_y]$ and all possible sequence of $x_n,y_n$ such that
\begin{itemize}
\item $x_n\in [px], y_n\in [py], x_n\neq p, y_n\neq p$, $x_n\to p$ or $y_n\to p$,
\item each pair $x_n,y_n$ are joined by a minimizing geodesic such that in the case $x_n\to p$ the curve $[x_ny_n]$
converges to a segment of $[py]$ or if $y_n\to p$ the curve $[x_ny_n]$ converges to a segment of $[px]$.
\end{itemize}
The lower limit taken on all sequences satisfying these conditions
\[
\measuredangle^{\text{ls}}[p^x_y] := \liminf \modang^0 (p^{x_n}_{y_n})
\]
is called the \textit{lower strong angle}.
Lower strong angles mainly serve to obtain symmetric estimates to those derived from the upper angle.

A \textit{triangle} $T$ in $X$ is the data of three points, called \textit{vertices}, and three minimizing geodesics, called \textit{edges}, joining them.
The notation $T^\kappa$ stands for a model triangle of $T$ in the model surface of constant curvature $\kappa$ and we write $|T^\kappa|$ its area measure.
Let $\delta_\kappa(T)$ denotes the \textit{relative excess} of $T$ with respect to $T^\kappa$, i.e.
\[
\delta_\kappa(T) := (\alpha + \beta + \gamma) - (\alpha_\kappa + \beta_\kappa + \gamma_\kappa)
\]
where $\alpha,\beta,\gamma$ are the angles at the vertices of $T$ and $\alpha_\kappa,\beta_\kappa,\gamma_\kappa$ are the angles of $T^\kappa$.
It will be also useful to express it as
\begin{equation}
\label{eq:RelativeExcessToDiff}
\delta_\kappa(T) = \delta_0(T) - \delta_0(T^\kappa).
\end{equation}
When the relative excess is computed with $\kappa=0$, we write simply it $\delta(T)$ and call it \textit{excess}.
When $\mathcal{T}$ is a finite family of non-overlapping triangles, then
\[
\delta(\mathcal{T}) := \sum_{T\in\mathcal{T}} \delta(T).
\]
With indeed $\delta(\mathcal{T})=0$ when $\mathcal{T}=\varnothing$.
A triangle $T$ is \textit{convex relative to the boundary} if no couple of points on the boundary of $T$ can be joined
by a curve lying outside $T$ and shorter than the part of the boundary joining the points.
A triangle is \textit{simple} if it is homeomorphic to a disc and convex relative to the boundary.

The next two propositions are general estimates for the distortion between angles in a metric space and model angle.
They first appeared in the paper \cite{alexandrowUberVerallgemeinerungRiemannschen1957} by Alexandrov.
They play a central role in the synthetic attempt to prove of the main theorems.

\begin{proposition}
\th\label{thm:AngleDistUPEstimate}
Let $T=[pqr]$ be a triangle in a metric space $(X,d)$.
Define
\[
\mu_\kappa := \sup_{\substack{x\in[pq], \\ y\in[pr]}} \delta_\kappa [pxy].
\]
Then
\[
\alpha - \alpha_\kappa \leq \mu_\kappa,
\]
where $\alpha$ is the upper angle at $p$ and $\alpha_\kappa$ the model angle.
\end{proposition}

\begin{proposition}
\th\label{thm:AngleDistLOWEstimate}
Let $T=[pqr]$ be a triangle in a metric space $(X,d)$ such that, for all $x\in[pq]$ and $y\in[pr]$, there exists a unique minimizing geodesic $[xy]$.
Assume strong angles exist between $[xy]$ and both $[xp], [yp]$.
For any $\kappa\in\R$, define
\[
\nu^{\text{ls}}_\kappa := \inf_{\substack{x\in[pq], \\ y\in[pr]}} \delta^{\text{ls}}_\kappa[pxy].
\]
where $\delta^{\text{ls}}_\kappa[pxy]$ is the relative excess computed with lower strong angles.
Then
\begin{equation}
\alpha^{\text{ls}} - \alpha_\kappa \geq \nu^{\text{ls}}_\kappa,
\end{equation}
where $\alpha^{\text{ls}}$ is the lower strong angle at $p$ and $\alpha_\kappa$ the model angle.
\end{proposition}

\subsection{Bounded integral curvature}
\label{subsec:BIC}

\begin{definition}
A metric space $S$ with distance function $d$ is called a \textit{surface with locally bounded integral curvature}, in short a \textit{BIC surface}, if
\begin{itemize}
\item $S$ is a connected, oriented, 2-dimensional smooth manifold without boundary,
\item $d:S\times S\to[0,+\infty)$ is an intrinsic distance inducing the manifold topology of $S$,
\item for every compact $K$ there exists a constant $C=C(K)>0$ such that $\delta(\mathcal{T})\leq C$
for any finite collection $\mathcal{T}$ of non-overlapping simple triangles contained in $K$.
\end{itemize}
\end{definition}

We fix a BIC surface $S$ for the rest of this section.
Let $\mu$ denote the 2-dimensional Hausdorff measure induced by $d$.
A signed measure of interest on these surfaces is the \textit{curvature measure} denoted $\omega$.
The positive and negative parts of this latter can be defined over open sets $U\subset S$ by
\[
\omega^+(U) := \sup_{\mathcal{T}} \sum_{T\in\mathcal{T}}\delta(T),\quad
\omega^-(U) := -\inf_{\mathcal{T}} \sum_{T\in\mathcal{T}}\delta(T),
\]
with both supremum taken on finite families $\mathcal{T}$ of non-overlapping simple triangles contained in $U$.
We extend $\omega^\pm$ to any other Borel subset $E\subset S$ by
\[
\omega^\pm(E) := \inf_{U\supset E~\text{open}} \omega^\pm(U),
\]
and so $\omega(E) = \omega^+(E) - \omega^-(E)$.
The measure $\omega^+$ is always locally finite, while local finiteness of $\omega^-$ is a deep fact of the theory
(see \cite[Theorem 15 p.134]{aleksandrovIntrinsicGeometrySurfaces2002}).
By \cite[Theorem 2.18]{rudinRealComplex1966}, $\omega$ is thus a signed Radon measure.
A point $p\in S$ is a \textit{cusp} if $\omega(\{p\})\geq 2\pi$.
This name comes the fact that, if $\omega(\{p\})> 2\pi$ then $p$ is a point at infinity and if $\omega(\{p\})=2\pi$
then $p$ may be at finite or infinite distance from other points.

For a simple closed curve $\gamma\subset S$ bounding a domain $D$ homeomorphic to a disc,
the \textit{Gauss-Bonnet formula} holds:
\[
\kappa_i(L) + \omega(D) = 2\pi,
\]
where $\kappa_i(L)$ is the turn of $L$ on the side of $D$ (see \cite[Chapter 13]{fillastreSubharmonic2023}).

These surfaces admit a local isothermal coordinates representation.
Let $D$ be a bounded domain of $\R^2$, let $\omega$ be a finite signed measure with compact support on $\R^2$ and let $h$ be a harmonic function
on $D$.
The pseudo-distance $d_{\omega, h}$, defined for every $x,y\in D$ by
\[
d_{\omega, h}(x,y) := \inf\left\{ \int_0^1 e^{u(\gamma(t))}|\gamma'(t)|\,dt\>\> \middle|\>\>
\gamma:[0,1]\to S~\text{rectifiable s.t.}~\gamma(0)=x, \gamma(1)=y \right\},
\]
where
\[
u(x) = -\frac{1}{\pi}\int_{\R^2} \ln|x-y|\dd\omega(y) + h(x),
\]
is called a \textit{$\delta$-subharmonic distance}.

\begin{theorem}
\th\label{thm:ExistenceSubharmonicCoor}
Let $(S,d)$ be a BIC surface and let $D\subset S$ be a domain whose closure is homeomorphic to a closed disc.
Suppose that $D$ is endowed with the intrinsic distance $\tilde d$ induced by the restriction of $d$ to $D$.
Then, there exists an isometry $\Phi: (D,\tilde d) \to (\Omega,\rho)$ such that $\rho$ is a $\delta$-subharmonic distance satisfying
$\rho=d_{\omega_D, h}$ where $\omega_D=\Phi_\#(\omega|_D)$ and $h$ is harmonic on $\Omega$.
\end{theorem}

\begin{proof}
See \cite[Theorem II Chap.7]{fillastreSubharmonic2023}.
\end{proof}

Reshetnyak also proved in \cite[Chap. 7]{fillastreSubharmonic2023} a convergence theorem for $\delta$-subharmonic distances.

\begin{theorem}
\th\label{thm:ReshDistCV}
Let $(\omega_n^1)$ and $(\omega_n^2)$ be two sequences of positive measures with supports contained in a same disc.
Suppose that $\omega_n^1\rightharpoonup\omega^1$ and $\omega_n^2\rightharpoonup\omega^2$.
Let $D\subset\C$ be the closure of a bounded domain whose boundary is the finite union of simple curves of bounded rotation.
Suppose that, for any $z\in D$,
\[
\omega^1(z) < 2\pi.
\]
Then, with $\omega:=\omega^1-\omega^2$ and $\omega_n:=\omega^1_n-\omega^2_n$, $(d_{\omega_n})$ converges uniformly to $d_\omega$ on $D$.
\end{theorem}

\subsection{Alexandrov spaces}

A full exposition of Alexandrov spaces can be found in \cite{alexanderAlexandrovGeometry2024}.

\subsubsection{Curvature bounded below}

\begin{definition}
A metric space $X$ is said to be $\CBB(\kappa)$ if for every quadruple of points $p,x_1,x_2,x_3\in X$
the following condition holds:
\[
\modang^\kappa(p^{x_1}_{x_2}) + \modang^\kappa(p^{x_2}_{x_3}) + \modang^\kappa(p^{x_3}_{x_1}) \leq 2\pi,
\]
or one of the angles $\modang^\kappa(p^{x_i}_{x_j})$ is not defined.
\end{definition}

\begin{proposition}
Let $\omega$ be a selective ultrafilter on $\N$.
Let $X_n$ be a $\CBB(\kappa_n)$ space for each $n$.
Assume that $X_n\to X$ and $\kappa_n\to\kappa$ as $n\to\omega$.
Then $X$ is $\CBB(\kappa)$.
In particular, the statement holds for Gromov-Hausdorff, uniform and Lipschitz convergence.
\end{proposition}

$X$ is said to be \textit{locally $\CBB(\kappa)$} if around every point $p\in X$ there is a neighborhood $U$
such that $U$ is $\CBB(\kappa)$.

When $X$ is a nice enough space it is possible to prove an equivalent statement to $\CBB(\kappa)$ condition.
This is the content of the next proposition.

\begin{proposition}
\th\label{alex:thm:CBBCV}
If $X$ is $\CBB(\kappa)$, then the following condition, called hinge comparison, holds:
For any hinge $[p_q^r]$, the angle $\measuredangle[p_q^r]$ is defined and
\[
\measuredangle[p_q^r] \geq \modang^\kappa(p_q^r).
\]
Moreover,
\[
\measuredangle[p_q^r] + \measuredangle[p_r^s] \leq \pi.
\]
Conversely, if $X$ is G-delta geodesic, then the hinge comparison implies $\CBB(\kappa)$.
\end{proposition}

CBB spaces also enjoy a local-to-global property.

\begin{theorem}
If $X$ is a complete length locally $\CBB(\kappa)$ space then it is a $\CBB(\kappa)$ space.
\end{theorem}

We say that $X$ is $\CBB(\kappa, n)$ if $X$ has Hausdorff dimension $n$ and $X$ is $\CBB(\kappa)$.

\begin{theorem}
For any $n\in\N$ and any $\kappa\in\R$, $\CBB(\kappa, n) \subseteq \RCD((n-1)\kappa, n)$.
Moreover, in the case of surfaces (i.e., $n=2$), the equality holds.
\end{theorem}

\begin{proof}
The paper \cite{petruninAlexandrovMeetsLottVillaniSturm2011} establishes that spaces satisfying $\CBB(\kappa, n)$
are contained in $\text{CD}((n-1)\kappa, n)$.

It remains to show the infinitesimal Hilbertianity.
We only give an outline of the proof as it is non-trivial and out of the scope of this work.
In \cite{kuwaeSobolevSpacesLaplacian2001}, the authors develop Sobolev spaces tailored to their work on $\CBB(\kappa, n)$ spaces.
They introduce a Dirichlet form $\mathcal{E}$, which is shown to be the Cheeger form in \cite[Section 2.2]{gigliHeatFlow2013}.
Since the Cheeger form is quadratic, it follows that $\CBB(\kappa, n)$ spaces are infinitesimally Hilbertian.

The reverse inclusion, when $n=2$, is proven in \cite{lytchakRicciCurvatureDimension2022}.
\end{proof}

\subsubsection{Curvature bounded above}

\begin{definition}
A metric space $X$ is said to be $\CAT(\kappa)$
if for every quadruple of points $p,q,x,y\in X$, one of the two following inequalities holds:
\[
\modang^\kappa(p_x^y) \leq \modang^\kappa(p^q_x) + \modang^\kappa(p_q^y)\quad\text{or}\quad
\modang^\kappa(q_x^y) \leq \modang^\kappa(q^p_x) + \modang^\kappa(q_p^y),
\]
or one of the corresponding model angles is not defined.
\end{definition}

$X$ is said to be \textit{locally $\CAT(\kappa)$} if around every point there is a neighborhood $U$
such that $U$ is $\CAT(\kappa)$.

\begin{proposition}
If $X$ is $\CAT(\kappa)$, then the following condition, called hinge comparison, holds:
For any hinge $[p_q^r]$, the angle $\measuredangle[p_q^r]$ is defined and
\[
\measuredangle[p_q^r] \geq \modang^\kappa(p_q^r).
\]
Conversely, if $X$ is $\varpi^\kappa$-geodesic, then the hinge comparison implies $\CAT(\kappa)$.
\end{proposition}

\section{Proof of \texorpdfstring{$\mathbf{\BICB\Rightarrow\CBB}$}{BICB => CBB}}

Let $(S,d)$ be a $\BICB(\kappa)$ surface with $d$ complete.
The proof employs a local-to-global strategy.
Thus we will first reduce the problem to a subharmonic distance on a domain using Reshetnyak's conformal charts
(see \thref{thm:ExistenceSubharmonicCoor}) and then apply CBB globalization on $S$ to obtain \thref{thm:BICBImpliesCBB}.

Let $\Omega\subset\R^2$ be an open disc such that the intrinsic distance is induced by the non-smooth metric $g=e^{2u}|\dd z|^2$
with $u$ a difference of two subharmonic functions.
Writing $\omega_\kappa:=\omega-\kappa\mu$, the hypotheses
\begin{equation}
\label{eq:lb}
\omega \geq \kappa\mu
\quad\Longleftrightarrow\quad \omega_\kappa\geq 0
\quad\Longleftrightarrow\quad -\Delta u \ge \kappa\,e^{2u}
\end{equation}
are equivalent.
Let first treat the hard case when $\kappa\leq0$.

For a standard mollifier $\rho_n\ge0$,
$\int\rho_n=1$, $\supp\rho_n\subset B_{1/n}$.
Set
\begin{equation}
\label{eq:fn}
f_n := \omega_\kappa * \rho_n \geq 0, \qquad f_n\in C^\infty(\R^2),\qquad f_n\,\dd x\dd y \rightharpoonup \omega_\kappa.
\end{equation}

Let $u_0:=u|_{\partial\Omega}$ (the boundary trace, finite, continuous and in $H^{1/2}(\partial\Omega)$
after slightly shrinking the radius of $\Omega$, see \thref{thm:admissibleBD}).
Define $u_n$ as the solution of the Dirichlet problem
\begin{equation}
\label{eq:liouville}
\begin{cases}
-\Delta u_n = \kappa\,e^{2u_n} + f_n &\text{in}~\Omega, \\
u_n = u_0 &\text{on}~\partial\Omega,
\end{cases}
\end{equation}
where $\Delta=\partial^2/\partial x^2 + \partial^2/\partial y^2$ is the usual Laplace operator and set $g_n:=e^{2u_n}|dz|^2$.
Equation~\eqref{eq:liouville} is exactly the prescription \textit{``the curvature measure of $g_n$ equals
$\kappa\,\vol_{g_n} +f_n\,\dd x\dd y$''}.

\begin{lemma}[Existence and lower bound]
\th\label{lem:LiouvilleSolExist}
Assume $\kappa\le 0$.
Then \eqref{eq:liouville} has a unique solution $u_n\in C^\infty(\Omega)\cap C(\overline\Omega)$, and
\begin{equation}
\label{eq:Klb}
\mathbb{K}_{g_n} = \kappa + e^{-2u_n} f_n \ge \kappa \qquad\text{pointwise in }\Omega,\ \text{for every }n.
\end{equation}
\end{lemma}

\begin{proof}
Existence and uniqueness: for $\kappa\leq0$ the equation \eqref{eq:liouville} is the Euler-Lagrange equation of the strictly convex, coercive
and weakly lower semicontinuous functional
\[
J(v) = \int_\Omega \left(\frac{1}{2}|\nabla v|^2 - \frac{\kappa}{2}e^{2v} - f_n\,v\right)\dd x\dd y, \qquad v\in u_0+H_0^1(\Omega),
\]
finite by Young and Moser-Trudinger inequalities \cite{cianchiMoserTrudinger2005}.
By \thref{thm:admissibleBD}, $u_0+H_0^1(\Omega)$ is a non-empty closed affine subset of $H^1(\Omega)$ made of $H^1$ functions with trace $u_0$.
Its unique minimiser $u_n$, existing by \cite[Corollary 3.23]{brezisFunctionalAnalysis2011} and unique by strict convexity, is the unique weak solution,
and elliptic regularity with bootstrapping give smoothness since $f_n\in C^\infty(\R^2)$.
Indeed, by Moser-Trudinger, the right-hand side of \eqref{eq:liouville} is in $L^p(\Omega)$, for any $p>1$, and by Calder\'on-Zygmund $u_n\in W^{2,p}(\Omega)$.
Sobolev embedding \cite[Theorem 7.26]{gilbargEllipticPDE2001} gives $u_n\in C^{1,\alpha}(\overline\Omega)$
and so $\kappa e^{2u_n} + f_n \in C^{1,\alpha}(\overline\Omega)\subset C^{0,\alpha}(\overline\Omega)$.
By \cite[Theorem 6.13]{gilbargEllipticPDE2001}, there is a unique solution
$\tilde u_n\in C(\overline\Omega)\cap C^{2,\alpha}(\Omega)$ of the Dirichlet problem $-\Delta\tilde u_n = \kappa e^{2u_n} + f_n$ on $\Omega$ and $u=u_n$ on $\partial\Omega$.
However, the function $w_n:=u_n-\tilde u$ is weakly harmonic on $\Omega$ and zero on the boundary.
So, by Weyl theorem, $w_n$ is smooth and $w_n\equiv 0$ everywhere on $\Omega$ by uniqueness.
Thus, $u_n\in C^2(\Omega)$.
By iterations of the Schauder interior estimate \cite[Theorem 6.17]{gilbargEllipticPDE2001}, $u_n\in C^\infty(\Omega)\cap C(\overline\Omega)$.

Identity \eqref{eq:Klb} is immediate from $K_{g_n}=e^{-2u_n}(-\Delta u_n)$ and \eqref{eq:liouville}, using $f_n\ge0$.
\end{proof}

\begin{lemma}[Mass below the threshold]
\th\label{lem:verify-bm}
There are $r_0>0$ and $\varepsilon_0>0$ such that
$f_n\big(B_{r_0}(z)\big)\le 2\pi-2\varepsilon_0$ for all $z\in\overline{\Omega}$ and all large $n$.
\end{lemma}

\begin{proof}
Regard $\omega_\kappa$ as a finite nonnegative Radon measure on $\R^2$ (extended by zero outside $\overline{\Omega}$).
Its atoms coincide with those of $\omega$, every atom has mass $<2\pi$ by assumption,
and for each $\eta>0$ only finitely many atoms have mass $\geq\eta$.
Define
\begin{equation}
\label{eq:Astrict}
A:=\sup_{x\in\R^2}\omega_\kappa(\{x\}) < 2\pi.
\end{equation}

\smallskip
We claim that,
\[
\lim_{r\to 0} \sup_{z\in \overline\Omega}\omega_\kappa\bigl(B_r(z)\bigr)
= \sup_{z\in \overline\Omega}\omega_\kappa(\{z\})
\leq A.
\]
Indeed, the map $r\mapsto\Phi(r):=\sup_{z\in \overline\Omega}\omega_\kappa(B_r(z))$ is nondecreasing, so
$L:=\lim_{r\to 0}\Phi(r)=\inf_{r>0}\Phi(r)$ exists.
For all $z\in\overline\Omega$ and $r>0$ we have $\omega_\kappa(B_r(z))\ge\omega_\kappa(\{z\})$, hence
$\Phi(r)\ge\sup_{z\in \overline\Omega}\omega_\kappa(\{z\})$ and therefore $L\ge\sup_{z\in \overline\Omega}\omega_\kappa(\{z\})$.

For the reverse inequality pick $r_k\searrow 0$ and $z_k\in\overline\Omega$ with
$\omega_\kappa\bigl(B_{r_k}(z_k)\bigr)\ge\Phi(r_k)-\tfrac1k\ge L-\tfrac1k$.
By compactness of $\overline\Omega$ we may assume $z_k\to z_\ast\in \overline\Omega$.
Fix $\rho>0$.
Since $r_k\to0$ and $z_k\to z_\ast$, we have $B_{r_k}(z_k)\subset B_\rho(z_\ast)$ for all large $k$,
so $\omega_\kappa\bigl(B_\rho(z_\ast)\bigr)\ge\omega_\kappa\bigl(B_{r_k}(z_k)\bigr)\ge L-\tfrac1k$.
Letting $k\to\infty$ gives $\omega_\kappa\bigl(B_\rho(z_\ast)\bigr)\ge L$.
Letting $\rho\to 0$ and using continuity from above of the finite measure $\omega_\kappa$
(note $\bigcap_{\rho>0}B_\rho(z_\ast)=\{z_\ast\}$) yields $\omega_\kappa(\{z_\ast\})\ge L$,
so $\sup_{z\in \overline\Omega}\omega_\kappa(\{z\})\ge L$.
Combining the two inequalities proves the equality, and $\sup_{z\in \overline\Omega}\omega_\kappa(\{z\})\le A$
is immediate from the definition of $A$.

By \eqref{eq:Astrict} we may set
\[
\varepsilon_0:=\frac{2\pi-A}{4}>0,
\qquad\text{so that}\qquad
A+2\varepsilon_0=2\pi-2\varepsilon_0.
\]
By the claim and \eqref{eq:Astrict} there is $r_1>0$ with
\[
\sup_{z\in \overline\Omega}\omega_\kappa\bigl(B_{r_1}(z)\bigr) \leq A+\varepsilon_0 = 2\pi-3\varepsilon_0.
\]
Set $r_0:=r_1/2$. Since $\supp\rho_n\subset B_{1/n}$, for every $z$
\[
\begin{aligned}
f_n\bigl(B_{r_0}(z)\bigr)
= \int_{B_{r_0}(z)}\left(\int_{\R^2}\rho_n(x-y)\dd\omega_\kappa(y)\right)\dd x
&= \int_{\R^2}\left(\int_{B_{r_0}(z)}\rho_n(x-y)\dd x\right)\dd\omega_\kappa(y) \\
&\leq \omega_\kappa\bigl(B_{r_0+1/n}(z)\bigr),
\end{aligned}
\]
because the inner integral never exceeds $\int_{\R^2}\rho_n=1$ and vanishes
unless $y\in B_{r_0+1/n}(z)$.
Hence, as soon as $1/n<r_1/2$ (so that $r_0+1/n<r_1$),
\[
f_n\bigl(B_{r_0}(z)\bigr)
\leq \omega_\kappa\bigl(B_{r_0+1/n}(z)\bigr)
\leq \sup_{w\in \overline\Omega}\omega_\kappa\bigl(B_{r_1}(w)\bigr)
\leq 2\pi-3\varepsilon_0
< 2\pi-2\varepsilon_0
\]
for every $z\in\overline\Omega$, where the second inequality uses
$B_{r_0+1/n}(z)\subset B_{r_1}(z)$ and $z\in\overline\Omega$.
This is the assertion.
\end{proof}

\begin{lemma}[Uniform global bound]
\th\label{lem:global}
For some $q>1$,
\begin{equation}\label{eq:globalLq}
\sup_n\big\|e^{2u_n}\big\|_{L^q(\Omega)}<\infty.
\end{equation}
In particular $\vol_{g_n}(\Omega)$ and $|\omega_n|(\Omega)$ are bounded uniformly in $n$.
\end{lemma}

\begin{proof}
For all $n$,
\begin{equation}\label{eq:upper}
u_n \leq P[u_0]+\mathcal{G}[f_n]\qquad\text{on }\Omega.
\end{equation}
Indeed, by the Green representation of \eqref{eq:liouville}, $u_n=P[u_0]+\mathcal{G}[\kappa e^{2u_n}+f_n]$.
Since $\kappa\le0$ and $\mathcal{G}\ge0$ the nonlinear term is nonpositive, giving \eqref{eq:upper}.

Thus, with $P[u_0]$ bounded on $\overline\Omega$, it suffices to bound $\sup_n\|e^{2\mathcal{G}[f_n]}\|_{L^q(\Omega)}$.
We estimate $e^{2\mathcal{G}[f_n]}$ on a half-ball about each point of $\overline\Omega$,
splitting the source by where its mass sits.

By Lemma~\ref{lem:verify-bm} there are $r_0,\varepsilon_0>0$ with $f_n(B_{r_0}(z))\le
2\pi-2\varepsilon_0$ for all $z\in\overline\Omega$ and all large $n$. Fix $z\in\overline\Omega$
and split $f_n=f_n'+f_n''$ with $f_n':=f_n\mathbf 1_{B_{r_0}(z)}$ and
$f_n'':=f_n\mathbf 1_{\Omega\setminus B_{r_0}(z)}$, so $\mathcal{G}[f_n]=\mathcal{G}[f_n']
+\mathcal{G}[f_n'']$.

The inner source has mass $\|f_n'\|_{L^1}=f_n(B_{r_0}(z))\le2\pi-2\varepsilon_0$, so the
Brezis-Merle inequality, see \cite{brezisUniform1991}, on $\Omega$ gives, for a fixed small $\delta$ with
$2q:=\frac{4\pi-\delta}{2\pi-2\varepsilon_0}>2$,
\[
\sup_n\int_\Omega e^{2q\mathcal{G}[f_n']} \leq \frac{4\pi^2}{\delta}(\operatorname{diam}\Omega)^2 .
\]
The outer potential $\mathcal{G}[f_n'']$ is harmonic on $B_{r_0}(z)$, and by the kernel
bound $0\le\mathcal{G}(\xi,\zeta)\le\frac1{2\pi}\log\frac{\operatorname{diam}\Omega}{|\xi-\zeta|}$,
for $\xi\in B_{r_0/2}(z)\cap\Omega$ and $\zeta\in\supp f_n''$ (so $|\xi-\zeta|\ge r_0/2$),
\[
0 \leq \mathcal{G}[f_n''](\xi)
\leq \frac1{2\pi}\log\frac{2\operatorname{diam}\Omega}{r_0}\, \|f_n''\|_{L^1} \leq C(r_0)\,M ,
\]
uniformly in $n$, with $M=\omega_\kappa(\Omega)$ and $\|f_n''\|_{L^1}\le\|f_n\|_{L^1}\le M$.
Hence on $B_{r_0/2}(z)\cap\Omega$,
\[
\int_{B_{r_0/2}(z)\cap\Omega}e^{2q\mathcal{G}[f_n]}
\leq e^{2qC(r_0)M}\int_\Omega e^{2q\mathcal{G}[f_n']}
\leq C,
\]
uniformly in $n$ and with $C$ independent of $z$.
The half-balls $\{B_{r_0/2}(z)\}_{z\in\overline\Omega}$ cover the compact $\overline\Omega$.
A finite subcover gives $\sup_n\|e^{2\mathcal{G}[f_n]}\|_{L^q(\Omega)}<\infty$, whence by
\eqref{eq:upper}
\[
\sup_n\int_\Omega e^{2qu_n}\ \le\ e^{2q\|P[u_0]\|_\infty}\,\sup_n\int_\Omega
e^{2q\mathcal{G}[f_n]}\ <\ \infty,
\]
which is \eqref{eq:globalLq}.

Finally $A_n(\Omega)\le|\Omega|^{1/q'}\|e^{2u_n}\|_{L^q(\Omega)}$
and $|\omega_n|(\Omega)\le\int_\Omega f_n+|\kappa|A_n(\Omega)\le M+|\kappa|\,C$ are bounded
uniformly in $n$.
\end{proof}

\begin{lemma}[$L^1$ convergence of the conformal factors]
\th\label{thm:u-conv}
$u_n\to u$ in $L^1_\loc(\Omega)$, and $e^{2u_n}\to e^{2u}$ in $L^1(\Omega)$.
\end{lemma}

\begin{proof}
\textit{Compactness.}
By \thref{lem:global}, $-\Delta u_n=\omega_n$ is bounded in $\mathcal{M}(\Omega)$.
The Stampacchia-Boccardo-Gallouët estimate, see \cite[Proposition 5.1]{ponceEllipticMeasure2016}, gives
$(u_n)$ bounded in $W^{1,p}(\Omega)$ for every $p<2$, hence precompact in $L^1(\Omega)$ by Rellich-Kondrachov.
Pass to a subsequence with $u_n\to u_*$ in $L^1(\Omega)$ and a.e.

\textit{Identification.}
Pass to the limit in \eqref{eq:liouville} in $\mathcal{D}'(\Omega)$:
$-\Delta u_n\xrightarrow[]{\mathcal{D}'}-\Delta u_*$, $\kappa e^{2u_n}\to\kappa e^{2u_*}$ in
$L^1_\loc$, and $f_n\rightharpoonup\omega_\kappa$, so
\[
-\Delta u_* = \kappa\,e^{2u_*} + \omega_\kappa,\qquad u_*=u_0\ \text{on}~\partial\Omega,
\]
We want to show that $u_*=u$.
Let $w:=u_*-u$.

\textit{Step 1 (exponential integrability).} We claim
\begin{equation}\label{eq:expLq}
e^{2u},\;e^{2u_*}\ \in\ L^q(\Omega)\qquad\text{for some }q>1 .
\end{equation}
Take the Riesz decomposition $u=\mathcal N[\omega^+]-\mathcal N[\omega^-]+h$ on a
disc containing $\overline\Omega$, where $\mathcal N[\mu](z)=\frac1{2\pi}\int\log
\frac1{|z-\zeta|}\dd\mu(\zeta)$, $\omega=\omega^+-\omega^-$ is the Jordan
decomposition, and $h$ is harmonic, hence bounded on $\overline\Omega$. Since
$\omega^-\ge0$,
\[
\mathcal N[\omega^-](z)\ \ge\ \tfrac1{2\pi}\Big(\log\tfrac1{\operatorname{diam}\Omega}\Big)\,
\omega^-(\Omega)\ >\ -\infty\qquad(z\in\Omega),
\]
so $e^{-2\mathcal N[\omega^-]}$ is bounded on $\Omega$.
Thus $e^{2u}\le C\,e^{2\mathcal N[\omega^+]}$ on $\Omega$.
By no-cusps assumption every point carries $\omega^+$-mass $<2\pi$, so by continuity of
$\omega^+$ from above each $z\in\overline\Omega$ admits a radius $r_z>0$ with
$\omega^+\big(\overline{B_{r_z}(z)}\big)<2\pi$.
The Brezis-Merle inequality (applied to the Newtonian potential of $\omega^+|_{B_{r_z}}$, the contribution of
$\omega^+|_{(\Omega\setminus B_{r_z})}$ being a bounded factor on $B_{r_z/2}$) gives $e^{2\mathcal N[\omega^+]}\in
L^{q_z}\big(B_{r_z/2}(z)\big)$ for some $q_z>1$.
The balls $B_{r_z/2}(z)$ cover the compact set $\overline\Omega$.
A finite subcover gives \eqref{eq:expLq} for $u$ with $q:=\min_i q_{z_i}>1$.
The same covering, applied to $u_n$ in place of $u$, is uniform in $n$: by \eqref{eq:globalLq} (Lemma~\ref{lem:global}),
$\sup_n\|e^{2u_n}\|_{L^q(\Omega)}<\infty$.
Since $e^{2u_n}\to e^{2u_*}$ a.e.\ along the subsequence, Vitali's theorem upgrades this to convergence in $L^q(\Omega)$.
In particular $F:=\kappa(e^{2u_*}-e^{2u})\in L^q(\Omega)$.

\textit{Step 2 ($w$ is the Green potential of an $L^q$ function).}
Each $u_n$ is smooth with $u_n|_{\partial\Omega}=u_0$.
By \thref{thm:admissibleBD} and potential regularity, $u\in W^{1,s}(\Omega)$ for $s<2$ with the same trace $u_0$.
Hence $w_n:=u_n-u\in W^{1,s}_0(\Omega)$, and subtracting $-\Delta u=\kappa e^{2u}+\omega_\kappa$ from \eqref{eq:liouville},
\[
-\Delta w_n=\underbrace{\kappa\big(e^{2u_n}-e^{2u}\big)}_{\in\,L^q(\Omega)}
+\underbrace{(f_n-\omega_\kappa)}_{\text{signed measure}}.
\]
The Laplace operator is linear, and a function in $W^{1,s}_0(\Omega)$ that is harmonic is smooth with zero trace, hence $\equiv0$.
So the $W^{1,s}_0$ solution is unique and equals the sum of the corresponding Green potentials,
\begin{equation}\label{eq:wn-split}
w_n=\mathcal{G}\!\big[\kappa(e^{2u_n}-e^{2u})\big]+\mathcal{G}[f_n-\omega_\kappa].
\end{equation}
Let $n\to\infty$ along the subsequence.
By the first step, $e^{2u_n}\to e^{2u_*}$ in $L^1(\Omega)$ with $\sup_n\|e^{2u_n}\|_{L^q(\Omega)}<\infty$.
Interpolation gives $\kappa(e^{2u_n}-e^{2u})\to F$ in $L^{q'}(\Omega)$ for some $q'\in(1,q)$, and since
$\mathcal{G}:L^{q'}(\Omega)\to C(\overline\Omega)$ is bounded with vanishing trace
(Lemma~\ref{lem:greenbound}), the first
term of \eqref{eq:wn-split} converges to $\mathcal{G}[F]$ in $C(\overline\Omega)$.
The second term tends to $0$ in $L^1(\Omega)$ by Stampacchia-Boccardo-Gallouët estimate.
As $w_n\to w$ in $L^1_\loc(\Omega)$, comparison of the limits yields
\begin{equation}
\label{eq:w-green}
w=\mathcal{G}[F]\ \in\ C(\overline\Omega),\qquad w|_{\partial\Omega}=0 .
\end{equation}
Here $\mathcal{G}[F]\in W^{2,q}(\Omega)\cap C(\overline\Omega)$ has zero trace, so the boundary
value in \eqref{eq:w-green} holds at every point of $\partial\Omega$.

\textit{Step 3 (sign and maximum principle).}
By \eqref{eq:w-green}, $-\Delta w=F$, and
\[
F\,w=\kappa\big(e^{2u_*}-e^{2u}\big)(u_*-u)\le0
\]
because $\kappa\le0$ and $t\mapsto e^{2t}$ is nondecreasing.
Thus $w$ is subharmonic on $\{w>0\}$ and superharmonic on $\{w<0\}$.
Being continuous on $\overline\Omega$ with $w|_{\partial\Omega}=0$, the maximum principle \cite[Theorem 2.3.1]{ransfordPotential1995}
on $\{w>0\}$ -- whose boundary lies in $\{w=0\}\cup\partial\Omega$ -- gives $\sup_\Omega w\le0$, and on
$\{w<0\}$ it gives $\inf_\Omega w\ge0$.
Hence $w\equiv0$, i.e.\ $u_*=u$.

As every subsequential limit equals $u$, the full sequence converges in $L^1_\loc(\Omega)$, and
$e^{2u_n}\to e^{2u}$ in $L^1(\Omega)$ by Step 1.
\end{proof}

\begin{lemma}[Convergence of the curvature measures]
\th\label{lem:measures}
$\dd A_n\to\dd A$ locally in total variation.
Moreover the curvature measures converge part by part: in the canonical nonnegative decomposition
\[
\omega_n = \omega_n^+ - \omega_n^-,\qquad
\omega_n^+ := f_n\,\dd x\dd y\ge0,\quad
\omega_n^- := |\kappa| \vol_{g_n}\ge0
\qquad(\kappa\le0),
\]
one has
\[
\omega_n^+\rightharpoonup\omega_\kappa,
\qquad
\omega_n^-\rightharpoonup|\kappa|\mu,
\]
with limit decomposition $\omega=\omega_\kappa-|\kappa|\mu$.
In particular $\omega_n\rightharpoonup\omega$.
\end{lemma}

\begin{proof}
$L^1$ convergence $e^{2u_n}\to e^{2u}$ (\thref{thm:u-conv}) yields $\vol_{g_n}\to\mu$ in total variation on compacts.
Testing against $C_b(\Omega)$ functions gives $\omega_n^-=|\kappa|\vol_{g_n}\rightharpoonup|\kappa|\,\mu$.
The positive part converges by construction: $\omega_n^+=f_n\dd x\dd y=\omega_\kappa*\rho_n\rightharpoonup \omega_\kappa$
by \eqref{eq:fn}.
Subtracting gives $\omega_n=\omega_n^+-\omega_n^-\rightharpoonup\omega_\kappa-|\kappa|\mu=\omega$, using
$\omega=\kappa\mu+\omega_\kappa$.
Both limit parts are nonnegative, and $\omega_\kappa(\{z\})=\omega^+(\{z\})<2\pi$ since $\kappa\mu$ is nonatomic.
\end{proof}

By Reshetnyak's \thref{thm:ReshDistCV}, the intrinsic distances $d_{g_n}\to d_g$ uniformly in $\Omega$.
Moreover, by \thref{alex:thm:CBBCV}, $\CBB(\kappa)$ is stable by uniform convergence.
At every point of $\Omega$ on a neighborhood small enough the distance $d$ coincide with $d_g$.
Thus, we deduce that $(S,d)$ is locally $\CBB(\kappa)$.
The proof is achieved using CBB globalization theorem.
(This is actually the only step where completeness is consumed.)

If $\kappa\ge0$, set $u_n=u*\rho_n$.
Then $-\Delta u_n=\omega*\rho_n$ and, by Jensen, $-\Delta u_n-\kappa e^{2u_n}=\omega_\kappa*\rho_n
+\kappa\,(e^{2u}*\rho_n-e^{2u_n})\ge0$, i.e.\ $K_{g_n}\ge\kappa$.
The area bound is then free: $\vol_{g_n}(\Omega)=\int e^{2(u*\rho_n)}\le\int e^{2u}*\rho_n\le\int e^{2u}$ by Jensen again.
All convergence statements above hold for this choice.
The arguments of the last paragraph hold too.
\thref{thm:BICBImpliesCBB} is now completely proven.

\subsection{Regularity and convergence lemmas}

Now we show that the initial claim that there is a ``nice'' domain such that the conformal factor is regular on the boundary.
For the next lemma, let $u$ be a $\delta$-subharmonic function on a domain $\Omega_0$.

\begin{definition}
A circle $C_r$ with $\overline{B_r(z_0)}\subset\Omega_0$ is called \textit{admissible} if
\begin{enumerate}[label=(\roman*)]
\item $|\omega|(C_r)=0$,
\item $C_r\cap\{p(\omega)=\pm\infty\}=\varnothing$,
\item $u\big|_{C_r}\in W^{1,3/2}(C_r)$.
\end{enumerate}
\end{definition}

\begin{lemma}[Admissible radius and boundary condition]
\th\label{thm:admissibleBD}
For almost every $r>0$ with $\overline{B_r(z_0)}\subset\Omega_0$ the circle $C_r(z_0)$ is admissible.
For such $r$, with $\Omega=B_r(z_0)$ and $u_0:=u\big|_{\partial\Omega}$,
\[
u_0~\text{is continuous on}~\partial\Omega, \qquad u_0\in H^{1/2}(\partial\Omega).
\]
In particular, the affine space $u_0+H^1_0(\Omega)$ is non-empty (by abuse of notation $u_0$ is an extesion
to $\Omega$ of $u_0$ on $\partial\Omega$).
\end{lemma}

\begin{proof}
\textit{1. Admissibility for a.e.\ $r$.}
The circles $\{C_r\}_{r>0}$ are pairwise disjoint.
On any compact annulus $\overline A\subset\Omega_0$ one has $|\omega|(\overline A)<\infty$, and a finite measure charges only
countably many members of a disjoint family.
Hence $\{r:|\omega|(C_r)>0\}$ is countable.
So (i) holds for all but countably many admissible radii.

The exceptional set $E:=\{p(\omega)=\pm\infty\}$ satisfies
$E\subseteq\{p(\omega^+)=-\infty\}\cup\{p(\omega^-)=-\infty\}$, a union of polar sets,
hence polar.
So $\dim_{\mathcal H} E=0$ and in particular $\mathcal H^1(E)=0$.
Moreover the radial projection $\pi(\zeta)=|\zeta-z_0|$ is $1$-Lipschitz, so $\pi(E)$ has Hausdorff
dimension $0$ and therefore Lebesgue measure zero in $(0,\infty)$.
Thus for a.e.\ $r$ one has $C_r\cap E=\varnothing$.

Fix $p=\tfrac{3}{2}$.
On an annulus $A=\{\rho_1<|z-z_0|<\rho_2\}\Subset\Omega_0$, polar coordinates give $ds=\rho\,d\theta$ on $C_\rho$
and $dz=\rho\,d\rho\,d\theta$, while the arc-length tangential derivative $\partial_\tau u=\tfrac1\rho\,
\partial_\theta u$ satisfies $|\partial_\tau u|\le|\nabla u|$.
Therefore
\[
\int_{\rho_1}^{\rho_2}\!\!\left(\int_{C_\rho} \bigl(|u|^{3/2}+|\partial_\tau u|^{3/2}\bigr)\,ds\right)d\rho
=\int_{A}\bigl(|u|^{3/2}+|\partial_\tau u|^{3/2}\bigr)\,dz \le\lVert u\rVert_{W^{1,3/2}(A)}^{3/2}<\infty .
\]
Hence the inner integral is finite for a.e.\ $\rho$.
Since $u\in W^{1,3/2}_\loc(\Omega_0)$ and the polar map $\Phi(\rho,\theta)=z_0+\rho e^{i\theta}$ is a diffeomorphism on $\overline A$,
the absolute continuity on lines theorem \cite[Theorem 10.35]{giovanniSobolev2009} applied to $u\circ\Phi$ in the angular variable gives,
for a.e.\ $\rho$, that $u|_{C_\rho}\in W^{1,3/2}(C_\rho)$ with $\tfrac{\dd}{\dd s}(u|_{C_\rho}) = (\partial_\tau u)|_{C_\rho}$.

The bad radii form a countable set united with two Lebesgue-null sets, hence a null set.
So a.e.\ $r$ with $\overline{B_r(z_0)}\subset\Omega_0$ is admissible.

\textit{2. Trace regularity.}
Fix an admissible $r$, parametrize $C_r$ by $\theta\mapsto z_0+re^{i\theta}$, and set $v(\theta):=u(z_0+re^{i\theta})$.
By~(iii), w.l.o.g.\ $v\in W^{1,3/2}(\mathbb{S}^1)$, i.e.\ $v,v'\in L^{3/2}(\mathbb{S}^1)$.
Writing $v(\theta)=\sum_{k\in\Z}\widehat v(k)e^{ik\theta}$ and applying the Hausdorff-Young inequality,
\[
\big(\widehat{v'}(k)\big)_k=\big(ik\,\widehat v(k)\big)_k \in\ell^3(\Z),
\qquad\text{i.e.}\qquad a_k:=|k|\,|\widehat v(k)|\in\ell^3(\Z) .
\]
By Hölder's inequality,
\[
\sum_{k\ne0}(1+|k|)\,|\widehat v(k)|^2
\leq 2\sum_{k\ne0}\frac{a_k^2}{|k|}
\leq 2\left(\sum_{k\ne0}a_k^3\right)^{2/3} \left(\sum_{k\ne0}|k|^{-3}\right)^{1/3}<\infty.
\]
With $|\widehat v(0)|<\infty$, and by \cite[Section 4.3]{taylorPDEI2011}, this shows that $v\in H^{1/2}(\mathbb{S}^1)$,
so $u_0:=u\big|_{C_r}\in H^{1/2}(\partial\Omega)$.
For the disk $\Omega=B_r(z_0)$, the trace operator $H^1(\Omega)\twoheadrightarrow H^{1/2}(\partial\Omega)$ is a surjection
\cite[Exercice 15.25]{giovanniSobolev2009},
hence the space $u_0+H^1_0(\Omega)$ is a non-empty affine subspace of $H^1(\Omega)$.

\textit{3. Continuity.}
Since $\dim C_r=1<\tfrac32=p$, Morrey's embedding on the one-dimensional circle gives
\[
W^{1,3/2}(C_r)\hookrightarrow C^{0,1/3}(C_r),
\]
so $u_0=u\big|_{C_r}$ admits a Hölder-continuous representative.
\end{proof}

\begin{lemma}[Green potentials of $L^{q'}$ data]
\label{lem:greenbound}
For every $q'>1$, with $q$ its conjugate exponent, the Green operator of $\Omega$ maps
$L^{q'}(\Omega)$ boundedly into $C(\overline{\Omega})$, the latter carrying the supremum
norm $\|f\|_{C(\overline{\Omega})}=\sup_{\overline{\Omega}}|f|$:
\[
\big\|\mathcal{G}[\phi]\big\|_{C(\overline{\Omega})}\ \le\ C_q\,\|\phi\|_{L^{q'}(\Omega)},
\qquad
C_q:=\sup_{z\in\overline\Omega}\big\|\mathcal{G}(z,\cdot)\big\|_{L^q(\Omega)}<\infty .
\]
Moreover $\mathcal{G}[\phi]\in C(\overline\Omega)$ with $\mathcal{G}[\phi]=0$ on $\partial\Omega$.
\end{lemma}

\begin{proof}
By Hölder's inequality, for every $z\in\overline\Omega$,
\[
\big|\mathcal{G}[\phi](z)\big|=\Big|\int_\Omega\mathcal{G}(z,\zeta)\,\phi(\zeta)\dd\zeta\Big|
\ \le\ \big\|\mathcal{G}(z,\cdot)\big\|_{L^q(\Omega)}\,\|\phi\|_{L^{q'}(\Omega)} .
\]
Write $D=\operatorname{diam}\Omega$, so $\Omega\subset B_D(z)$ and $0\le\mathcal{G}(z,\zeta)\le
\frac1{2\pi}\log\frac{D}{|z-\zeta|}$.
In polar coordinates centered at $z$, with the substitution $s=\log(D/r)$,
\[
\int_\Omega\left(\log\frac{D}{|z-\zeta|}\right)^{q}\dd\zeta
\leq 2\pi\int_0^D\left(\log\frac{D}{r}\right)^{q}r\dd r
= 2\pi D^2\!\int_0^\infty s^{q}e^{-2s}\dd s
= \frac{2\pi D^2\,\Gamma(q+1)}{2^{\,q+1}},
\]
independent of $z$.
Hence $C_q^{\,q}\le(2\pi)^{1-q}2^{-(q+1)}D^2\Gamma(q+1)<\infty$, and the asserted bound follows.
The logarithmic singularity lies in $L^q(\Omega)$ for
every $q<\infty$, so the bound holds for all $q'>1$ with no further restriction.

For continuity, fix $z_0\in\overline\Omega$ and let $z\to z_0$. Then $\mathcal{G}(z,\zeta)
\to\mathcal{G}(z_0,\zeta)$ for $\zeta\ne z_0$.
And the family $\{\mathcal{G}(z,\cdot)^q\}_{z\in \overline\Omega}$ is uniformly integrable, since the bound above has a radially
decreasing majorant about $z$, whence $\sup_z\int_E\mathcal{G}(z,\cdot)^q\dd\zeta\le
2\pi\int_0^{(|E|/\pi)^{1/2}}(\log\frac{D}{r})^q r\dd r\to0$ as $|E|\to0$.
By Vitali, $\mathcal{G}(z,\cdot)\to\mathcal{G}(z_0,\cdot)$ in $L^q(\Omega)$, so $\mathcal{G}[\phi](z)\to
\mathcal{G}[\phi](z_0)$ by Hölder.
Thus $\mathcal{G}[\phi]\in C(\overline\Omega)$.
Finally $\mathcal{G}(\xi,\cdot)\equiv0$ for $\xi\in\partial\Omega$, so $\mathcal{G}[\phi]|_{\partial\Omega} =0$.
\end{proof}

\section{Attempt with synthetic geometry and its limitations}

In this section we try to approach the problems $\BICB\Rightarrow\CBB$ and $\BICA\Rightarrow\CAT$ using only metric geometry.
This technique is limited by the fact that upper angles -- which are the one used by $\CBB/\CAT$ theory -- only see the metric structure.
On the other side, the Gauss-Bonnet theorem is of a topological nature and relates the \emph{sector angles} to the curvature of the triangle.
But sector angles depend on a choice of side of the hinge that the upper angle cannot catch
(see \cite[Theorem 6 p.120]{aleksandrovIntrinsicGeometrySurfaces2002}).

We were not able to pursue further in this direction however the techniques remain worth of interest.
In the case of $\BICB$, it is possible to show with metric method that $\CBB$ condition holds for convex triangles and that they are non-branching.
For $\BICA$ surfaces, it works when $\kappa<0$.

\subsection{Area and curvature}

To ensure that our analysis is consistent with the monograph \cite{aleksandrovIntrinsicGeometrySurfaces2002} of Alexandrov and Zalgaller,
we show that their definition of area measure corresponds to the 2-dimensional Hausdorff measure.
Let $S$ be a BIC surface.
In \cite{aleksandrovIntrinsicGeometrySurfaces2002}, the authors first define a pre-measure $\sigma_0$ on the set of polygons $P$ of $S$ by
\[
\sigma_0(P) := \lim_{\diam\mathcal{T}\to 0} \sum_{T\in\mathcal{T}} |T^0|
\]
where $\mathcal{T}$ are triangulations of $P$.
The \textit{area} $\sigma$ is then defined as an outer-measure by
\[
\sigma(E) := \inf_{O\supseteq E}\sup_{P\subseteq O} \sigma_0(P),\quad E\subseteq S.
\]
The next proposition gather elementary properties of $\sigma$.

\begin{proposition}
\th\label{thm:ArProp}
The following points hold:
\begin{enumerate}[label=(\roman*)]
\item
The limit in $\sigma_0(P)$ exists and doesn't depend on the chosen coverings of $P$.
\item
$\sigma$ is a measure on the Borel sets.
\item
Points, shortest curves and geodesics are null sets for $\sigma$.
\item
If $P$ is a polygon, then $\sigma(P)=\sigma_0(P)$.
\item
If $E\subseteq S$ has an interior point, then $\sigma(E)>0$.
\end{enumerate}
\end{proposition}

\begin{proof}
See \cite[Theorem 2 and 4-7 Chap. VIII]{aleksandrovIntrinsicGeometrySurfaces2002}.
\end{proof}

\begin{lemma}
\th\label{thm:ArMeasHausdorff}
Let $(S,d)$ be BIC surface.
The area measure $\sigma$ coincides with the 2-dimensional Hausdorff measure on Borel sets.
\end{lemma}

This result was already proved by M. B. Stratilatova in \cite{stratilatovaArea1957}.
Since the original article could not be located, we present an new proof here.
It is an improvement of \cite[Proposition 1.3]{chowdhuryCATKSurfaces2025} as
we do not make any assumption on curvature measure and presence of cusps.
This also answer \cite[Remark A.14]{thomasrichardCanonicalSmoothing2018}.

\begin{proof}
First, let admit that $(S,d)$ is a polyhedral surface and let $\mathcal{Z}$ be a triangulation of $S$
by triangles each one isometric to a plane triangle.

Let $P$ be a polygon in $S$.
We can triangulate $P$ such that each triangle is contained in a triangle of $\mathcal{Z}$.
However, on Euclidean triangles, $\sigma_0=\mathcal{H}^2=\mathcal{L}^2$.
Thus, by gluing parts together, we obtain that $\sigma_0(P)=\mathcal{H}^2(P)$ for every polygon $P$.

Let $O$ be an open subset of $S$.
For any $Z\in\mathcal{Z}$, trivially
\[
\sigma(O\cap Z^\circ)
= \sup_{P\subset O\cap Z^\circ}\sigma_0(P)
\leq \sup_{K\subset O\cap Z^\circ} \mathcal{H}^2(K)
= \mathcal{H}^2(O\cap Z^\circ)
\]
because every polygon is a compact and $\mathcal{L}^2$ is a Radon measure.
To conclude, we need to show the reverse inquality.
For that, we show that for every $K\subset O\cap Z^\circ$, there is a polygon $P$ containing $K$.
Since $K$ is compact and the complement of $O\cap Z^\circ$ in $\R^2$ is closed, then
\[
\varepsilon := d(K, (O\cap Z^\circ)^\complement) > 0.
\]
We cover the plane by a grid of squares of side length $\ell<\varepsilon/\sqrt{2}$.
Let $P$ be the union of all squares that intersect $K$.
$P$ is a polygon as its boundary is made of segments for the Euclidean metric.
It is clear that $K\subset P$ as every point of $K$ is contained in a grid square that intersect $K$.
Moreover, $P\subset Z^\circ$.
Indeed, let $x$ be any point in $P$.
Then $x$ belongs to some grid square that meets $K$ and so there exists $y\in K$ in the same square.
Since the diameter of the square is $\ell\sqrt{2} < \varepsilon$,
\[
d(x,K) \leq d(x,y) \leq \ell\sqrt{2} < \varepsilon = d(K, (O\cap Z^\circ)^\complement),
\]
so $x\not\in (O\cap Z^\circ)^\complement$.
This proves that $\sigma=\mathcal{H}^2$ on polyhedral surfaces.

Now we suppose that $(S,d)$ is a compact BIC surface without cusps.
By \cite[Lemma 6]{buragoBiLipschitz2005}, there exists a sequence $(d_i)$ of polyhedral metrics
converging to $d$ in the Lipschitz sens and such that $\omega_{d_i}^\pm\rightharpoonup\omega^\pm_d$.
The sequence $(d_i)$ is uniformly convergent and, by weak convergence, we get that $|\omega_{d_i}|(S)\to|\omega_d|(S)$
so $(|\omega_{d_i}|(S))$ is uniformly bounded.
Thus, we apply \cite[Theorem 9 p.269]{aleksandrovIntrinsicGeometrySurfaces2002} to obtain $\sigma_{d_i}\rightharpoonup\sigma_d$.
For any $E\subseteq S$ with $\partial E=\emptyset$,
\[
(1-\varepsilon)^2 \mathcal{H}^2_d(E) \leq \mathcal{H}^2_{d_i}(E) \leq (1+\varepsilon)^2 \mathcal{H}^2_d(E)
\]
Letting $\varepsilon\to 0$ this shows $\mathcal{H}^2_{d_i}(E)\to\mathcal{H}^2(E)$ and so $\mathcal{H}^2_{d_i}\rightharpoonup\mathcal{H}^2_d$.
However, in the last paragraph, we proved that $\mathcal{H}^2_{d_i}=\sigma_{d_i}$.
By uniqueness of the weak limit, $\mathcal{H}^2_d=\sigma_d$.

For the last case, we suppose that $(S,d)$ is an arbitrary BIC surface.
By \cite[Theorem 1 p.58]{aleksandrovIntrinsicGeometrySurfaces2002}, one can find at every point $p$ in $S$ a geodesically convex polygonal neighborhood $P$ homeomorphic to a disc
with diameter lower than an arbitrary constant.
Supposing that $P$ does not contain any cusp,
the space $(H\# P, \tilde d)$ is a compact BIC surfaces without cusps made of the gluing along the boundaries of $P$
and a hemisphere $H\subset\mathbb{S}^2$ with same radius as the boundary of $P$.
Thus, on $H\# P$, $\mathcal{H}^2_{\tilde d}=\sigma_{\tilde d}$ and so by isometry $\mathcal{H}^2_d|_P = \sigma_d|_P$.

Let $K\subseteq S$ be a compact set.
By compactness, $K$ can contains only a finite number of cusps $c_1,\dots, c_k$.
We build an exhaustion sequence $(K_n)$ of $K\setminus\{c_j\}$ by defining
\[
K_n := K \setminus \bigcup_{i=1}^k B_d(c_i, 1/n),\quad n\geq 1.
\]
The $K_n$ are compact and so we can cover them by finitely many polygonal neighborhood $(P_q)_{1\leq q\leq m}$
small enough such that any of them do not contain any cusp.
We obtain
\begin{align*}
\mathcal{H}^2_d(K_n) &= \sum_{q=1}^m \mathcal{H}^2_d(P_q\cap K_n) - \sum_{1\leq q<r\leq m} \mathcal{H}^2_d(P_q\cap P_r \cap K_n) \\
&= \sum_{q=1}^m \sigma_d(P_q\cap K_n) - \sum_{1\leq q<r\leq m} \sigma_d(P_q\cap P_r \cap K_n) \\
&= \sigma_d(K_n).
\end{align*}
We have constructed an increasing sequence $ K_1 \subseteq  K_2 \subseteq \dots \subseteq K\setminus\{c_j\}$ such that
$K\setminus\{c_j\}=\cup_i K_n$ and $\mathcal{H}^2_d(K_n)=\sigma_d(K_n)$.
Thus
\[
\mathcal{H}^2_d(K\setminus\{c_j\}) = \lim_{n\to\infty} \mathcal{H}^2_d(K_n) = \lim_{n\to\infty} \sigma_d(K_n) = \sigma_d(K\setminus\{c_j\}).
\]
However, because singletons have zero measure for $\mathcal{H}^2_d$ and $\sigma_d$,
\[
\mathcal{H}^2_d(K) = \mathcal{H}^2_d(K\setminus\{c_j\}) = \sigma_d(K\setminus\{c_j\}) = \sigma_d(K).
\]
Thus $\mathcal{H}^2_d=\sigma_d$ on compact subsets of $S$.
This achieves the proof.
\end{proof}

The next proposition is a key tool in the proof of Theorems \ref{thm:BICBImpliesCBB} and \ref{thm:BICAImpliesCAT}.
It relates the distortion of the area of a triangle to the curvature.
This result is of an independant interest in surface theory.
Following Alexandrov and Zalgaller, we adopt the notation $T_-$ for the interior of a triangle $T$ and
\[
\tilde\omega(T) := \omega(T_-) + \sum_{i=1}^3 \tau_i^-,
\]
where the $\tau_i$ are the turns of the edges of $T$ toward the interior.

\begin{proposition}
Let $T$ be a triangle in $S$ that is homeomorphic to a disc with $\diam(T)<\ell$.
Then,
\begin{equation}
\label{eq:AreaDistortion}
- \frac{1}{2}\tilde\omega^-(T)\ell^2 \leq \mu(T) - |T^0| \leq \frac{1}{2} \omega^+(T_-) \ell^2.
\end{equation}
\end{proposition}

\begin{proof}
From \cite[Lemma 4 p.260]{aleksandrovIntrinsicGeometrySurfaces2002}, we have that, for any $\varepsilon$,
there exists a arbitrarily fine triangulation of $T$ such that
\[
-\frac{1}{2}\tilde\omega^-(T)\ell^2-\varepsilon \leq |Q| - |T^0| \leq \frac{1}{2}\omega^+(T_-)\ell^2 + \varepsilon
\]
where $Q$ is a development in the plane of the triangulation and $|Q|$ its area measure.
Then, letting $\varepsilon\to 0$, $|Q|\to\sigma_0(T)=\sigma(T)=\mu(T)$ by \thref{thm:ArMeasHausdorff} and \thref{thm:ArProp}.
\end{proof}

\subsection{Curvature bounded below}

The next proposition rules out cusps from surfaces with non-negative curvature.
This result was already known by Reshetnyak (see \cite[p.130]{reshetnyakGeometryIV1993}).
We recall it for completeness.

\begin{lemma}
\label{lm:NotTooCuspsy}
On a $\BICB(0)$ surface $S$, $\omega(\{p\})<2\pi$ for every $p\in S$.
\end{lemma}

\begin{proof}
Recall that $p\in S$ is said to be \textit{at infinity} if for every $q\in S, d(q,p)=+\infty$.
From \cite[Theorem 4.110]{fillastreSubharmonic2023}, we know that a point with $\omega(\{p\})>2\pi$ is at infinity
which is not possible by definition.
Local finiteness of $\omega$ implies that the set of cusps is discrete.
Being a point at infinity is a local property.
We can restrict our analysis to a closed disc $S$ of the plane with $\diam(S)<1$ for the Euclidean distance and containing
only the origin as cusps.

The function $h$ is harmonic on $S$, so there is a constant $C>0$ such that,
\[
\begin{aligned}
\exp(u(z)) &= \exp\left( -\frac{1}{2\pi} \int_S \ln|z-\zeta|\,d\omega(\zeta) + h(z) \right) \\
&= |z|^{-1} \exp\left( -\frac{1}{2\pi} \int_{S\setminus\{0\}} \ln|z-\zeta|\,d\omega(\zeta) + h(z) \right) \\
&= C|z|^{-1}.
\end{aligned}
\]
Hence every curve starting at the origin will have an infinite length.
This is contradictory because the distance of a BIC surface is supposed to be finite.
\end{proof}

\begin{theorem}
\th\label{thm:CmBBDiscCvxImpliesHinge}
Suppose that $S$ is a $\BICB(\kappa)$ surface for some $\kappa\in\mathbb{R}$.
For each $\kappa'<\kappa$, there exists a $\ell>0$ such that for any triangle $T=[pqr]$ of $S$ satisfying:
\begin{itemize}
\item $T$ is homeomorphic to a disc,
\item $\diam(T)<\ell$,
\item the sectors inside the triangle at the vertices of $T$ are convex,
\end{itemize}
then the upper angle at $p$ is not less than the model angle for $\kappa'$, i.e.
\[
\measuredangle [p^q_r] \geq \modang^{\kappa'} (p^q_r).
\]
\end{theorem}

The following proof is inspired by the proof of \cite[Lemma 1 p.379]{aleksandrovIntrinsicGeometryConvex2006}.

\begin{proof}
Let $\omega$ denote the curvature measure of $S$.

By a procedure of excision and pasting from the proof of \cite[Theorem 10 p.216]{aleksandrovIntrinsicGeometrySurfaces2002},
we can admit that every triangle $[pxy]$ are triangles homeomorphic to a disc and absolutely convex.
We recall it here for completeness.

We excise $T$ from the surface $S$.
If $T$ has an interior tail, we cut the tail until the splitting point.
If $T$ has an exterior tail, we past along the side containing the tail the excision of an arbitrarily narrow two-gon
formed by two great circles on the sphere.
We now glue along the boundary of $T$ a semi-infinite right cylinder with the same perimeter as $T$.
After this excision-pasting procedure, in the resulting manifold endowed with the induced intrinsic distance,
the new triangle $\tilde T$ is homeomorphic to a disc and absolutely convex.
This procedure does not change the angle at $p$ as the sector at $p$ was convex by hypothesis.
This also does not change the model angles at $p$ as the geodesics between the vertices remain the sames.
The curvature measure and the area measure on the interior of $\tilde T$ stay also unchanged.
The conditions of \cite[Theorem 9 p.41]{aleksandrovIntrinsicGeometrySurfaces2002} are fullfilled for the excised triangle.

As $S$ is a BIC surface, angles coincide with lower strong angles.
So for the clarity of the proof we can simply work with angles.

Let $\kappa$ and $\kappa'$ be as prescribed in the statement of the theorem.
One can observe in the inequality
\[
\alpha - \alpha_\kappa \geq \nu_\kappa
\]
of \thref{thm:AngleDistLOWEstimate} that the proof is achieved if $\nu_{\kappa'}$ is shown to be non-negative.
In other words, if $\delta_{\kappa'}[pxy]\geq 0$ for every $x\in[pq], y\in[pr]$.
Having \eqref{eq:RelativeExcessToDiff} and Gauss-Bonnet theorem in mind, it boils down to showing that, for every $x\in[pq], y\in[pr]$,
\[
\delta_0[pxy] \geq \kappa' \left|[pxy]^{\kappa'}\right|.
\]
Starting with the curvature hypothesis $\omega\geq\kappa\mu$ applied to any triangle $T$,
\begin{equation}
\label{eq:HypothesisCurvBelow}
\delta_0(T)=\omega(T_-)\geq\kappa\mu(T)
\end{equation}
we will estimate from below the right-hand term $\kappa\mu(T)$.
The rest of the proof splits in two parts, one for the case $\kappa$ non-negative and one for $\kappa$ negative.

Recall that $\tilde\omega^-(T)=\omega^-(T_-)$ because the turn of shortest arcs is zero on $\BICB$ surfaces.

$(\kappa\geq 0)$
Directly from the area distortion estimate \eqref{eq:AreaDistortion} applied once on $T$ and once on $T^\kappa$, we have
\[
\begin{aligned}
\mu(T) &\geq |T^0| - \frac{1}{2}\omega^-(T)\ell^2 \\
&\geq |T^\kappa| - \frac{1}{2}\omega^+(T^\kappa)\ell^2 \\
&= \left( 1 - \frac{1}{2}\kappa \ell^2 \right) |T^\kappa|.
\end{aligned}
\]
When we plug it back into \eqref{eq:HypothesisCurvBelow}, it gives us
\[
\delta_0(T) \geq \kappa \left( 1 - \frac{1}{2}\kappa \ell^2 \right) |T^\kappa|.
\]
So by choosing $d$ such that
\[
\kappa' = \kappa \left( 1 - \frac{1}{2}\kappa \ell^2 \right)
\]
and by the fact that $|T^\kappa|\geq|T^{\kappa'}|$, we get
\[
\delta_0(T) \geq \kappa' |T^{\kappa'}|
\]
which is the desired inequality.

$(\kappa<0)$
The reasoning mirrors that of the non-negative case.
Again, by area distortion estimate \eqref{eq:AreaDistortion}, we have
\[
\begin{aligned}
\mu(T) &\leq \frac{1}{2}\omega^+(T_-)\ell^2 + |T^0| \\
&\leq \frac{1}{2}(\delta_0(T)+\omega^-(T_-))\ell^2 + \frac{1}{2} \omega^-(T^{\kappa'})\ell^2 + |T^{\kappa'}| \\
&\leq \frac{1}{2}(\delta_0(T)+|\kappa|\mu(T))\ell^2 + |T^{\kappa'}| \left( 1 + \frac{1}{2}|\kappa'|\ell^2 \right).
\end{aligned}
\]
Thus
\[
\left( 1 - \frac{1}{2} |\kappa|\ell^2 \right) \mu(T) \leq
\frac{1}{2}\delta_0(T)\ell^2 + |T^{\kappa'}| \left( 1 + \frac{1}{2}|\kappa'|\ell^2 \right)
\]
And so, for $\ell$ small enough, we plug it back into \eqref{eq:HypothesisCurvBelow},
\[
\delta_0(T) \geq \frac{1}{1 - \frac{1}{2} |\kappa|\ell^2}
\left( \frac{1}{2}\kappa\delta_0(T)\ell^2 + \kappa|T^{\kappa'}| \left( 1 + \frac{1}{2}|\kappa'|\ell^2 \right) \right).
\]
Simplifying the last inequality,
\[
\delta_0(T) \geq \kappa \left( 1 + \frac{1}{2} |\kappa'| \ell^2 \right) |T^{\kappa'}|.
\]
By shrinking $\ell$ even more if necessary, we let
\[
\kappa'' := \kappa \left( 1 + \frac{1}{2} |\kappa'| \ell^2 \right) > \kappa'
\]
to obtain
\[
\delta_0(T) \geq \kappa'' |T^{\kappa'}| \geq \kappa'|T^{\kappa'}|
\]
which is the desired inequality.
This achieves the proof.
\end{proof}

\begin{lemma}
\th\label{lm:NonBranching}
Assume $S$ is a $\BICB(\kappa)$ surfaces.
Then $S$ is non-branching.
\end{lemma}

\begin{proof}
Let $\omega$ denote the curvature measure of $S$.
Let $\gamma_0, \gamma_1$ be two branching geodesics starting at a point $p$ and ending at the points $x,y$ respectively.
We define
\[
t_{\text{max}} := \sup\{ t | \gamma_0(s) = \gamma_1(s)~\text{for all}~s\in[0,t)\}.
\]
By continuity, $\gamma_0(t_\text{max})=\gamma_1(t_\text{max})$.
Let's call this splitting point $z$.
We have, by definition of the turns,
\[
\kappa_l(\gamma_i) + \kappa_r(\gamma_i) = \omega(\gamma_i^\circ)
\]
where $\gamma_i^\circ$ stands for the image of $\gamma_i$ without the endpoints.
However, by no cusp assumption or lemma \ref{lm:NotTooCuspsy} and \cite[Theorem 4.121]{fillastreSubharmonic2023},
the turns are negatives so
\[
\kappa_l(\gamma_i)=\kappa_r(\gamma_i)=\omega(\gamma_i^\circ)=\omega(\{z\})=0.
\]
This shows that the complete angle at $z$ is $2\pi$.
The formula given by \cite[Theorem 4.156]{fillastreSubharmonic2023} applied on both $\gamma_i$ gives us that the angles at
$z$ between the two branches of $\gamma_i$ cut at $z$ is $\pi$.
Therefore the angle of the hinge $[z^x_y]$ must be $0$.
The sector induced by the hinge $[z^x_y]$ is geodesically convex.
Indeed, if a geodesic between two points in the sector crosses one of the $\gamma_i$ then it is possible to cut it with a segment
of $\gamma_i$ to shorten it which is contradictory.
Now we consider two points $q,r$ respectively on each of the branches of the hinge.
Two cases can happen.
First, the geodesic between $q$ and $r$ is fully included in $[zx]\cup[zy]$.
This is impossible otherwise, again using \cite[Theorem 4.156]{fillastreSubharmonic2023},
the angle of the sector will be $\pi$ which is impossible.
On the other hand, if we take $q,r$ close enough to $z$ it is possible to get a triangle homeomorphic to a disc
with a convex sector at $x$.
But by hinge comparison we have,
\[
\measuredangle[z^x_y] = 0 \Rightarrow
\modang^\kappa\left(z^{\operatorname{geod}_{[zx]}(t)}_{\operatorname{geod}_{[zy]}(t)} \right) = 0
\]
for $t$ small enough.
This means that $\gamma_0$ and $\gamma_1$ coincides after $z$, which contradicts the hypothesis.
\end{proof}

\subsection{Curvature bounded above}

\begin{theorem}
\th\label{thm:CmBADiscImpliesHinge}
Suppose that $S$ is a $\BICA(\kappa)$ surface for some $\kappa \leq 0$.
Then, for each $\kappa'>\kappa$, there exists a $\ell>0$ such that for any triangle $T=[pqr]$ of $S$ homeomorphic to a disc
and whose diameter is less than $\ell$, the upper angle is not more than the model angle for $\kappa'$, i.e.
\[
\measuredangle [p^q_r] \leq \modang^{\kappa'} (p^q_r).
\]
\end{theorem}

\begin{proof}
The proof follows similar lines as for \thref{thm:CmBBDiscCvxImpliesHinge} using opposite inequalities
and \thref{thm:AngleDistUPEstimate}.
For the sake of completeness we write it entirely.

The curvature measure on the interior of the excise triangle glued with a 2-gon from the sphere will satisfies the curvature bound
with an arbitrary small positive change appearing.
However, the steps remains unchanged.

Let $[pqr]$ be a triangle in $S$.
\thref{thm:AngleDistUPEstimate} shows that if
\[
\delta_0[pxy] \leq \kappa' \left|[pxy]^{\kappa'}\right|
\]
holds for any $x\in[pq], y\in[pr]$, then the conclusion of the theorem follows.
Starting with the curvature hypothesis $\omega\leq\kappa\mu$ applied to any triangle $T$ homeomorphic to a disc,
\begin{equation}
\label{eq:HypothesisCurvAbove}
\delta_0(T)\leq\omega(T_-)\leq\kappa\mu(T)
\end{equation}
we will estimate from above the right-hand term $\kappa\mu(T)$.
The rest of the proof splits in two parts, one for the case $\kappa=0$ and one for $\kappa$ negative.

$\boldsymbol{(\kappa = 0)}$
This case is straightforward.
By curvature hypothesis, $\delta_0(T)\leq 0$ and so, as $\kappa'$ is positive,
\[
\delta_0(T) \leq \kappa'|T^{\kappa'}|.
\]

$\boldsymbol{(\kappa < 0)}$
Without loss of generality, we can assume that $\kappa<\kappa'<0$.
We have that
\begin{align*}
\mu(T) &\geq -\frac{1}{2}\tilde\omega^-(T)\ell^2 + |T^0| \\
&\geq -\frac{1}{2}\tilde\omega^-(T)\ell^2 + |T^{\kappa'}| \\
&= \frac{1}{2}\tilde\delta_0(T) \ell^2 + |T^{\kappa'}| \\
&\geq \frac{1}{2}\delta_0(T) \ell^2 + |T^{\kappa'}|.
\end{align*}
Then
\[
\delta_0(T) \leq \kappa \left( \frac{1}{2}\delta_0(T)\ell^2 + |T^{\kappa'}| \right)
\]
and so
\[
\delta_0(T) \leq \frac{\kappa}{1+\frac{1}{2}|\kappa| \ell^2} |T^{\kappa'}|
\]
Thus, by taking $\ell$ small enough, we obtain
\[
\delta_0(T) \leq \kappa' |T^{\kappa'}|.
\]
This achieves the proof.
\end{proof}

\begin{lemma}
\th\label{thm:CmBBunicitygeod}
Let $S$ be a $\BICA(\kappa)$ surface with $\kappa<0$.
Then, for any points $p,q\in S$ contained in a domain homeomorphic to a disc, there is a unique geodesic between the two points.
\end{lemma}

\begin{proof}
Let $\omega$ denote the curvature measure of $S$.
Let $p,q\in S$ be two points joined by two different geodesics $\gamma_1, \gamma_2$ contained in a domain homeomorphic to a disc.
As these two geodesics are different, together they form a Jordan curve enclosing a domain $D$.
Let call this enclosing curve $\gamma$.
Then by Gauss-Bonnet theorem
\[
\tau_l(\gamma) + \omega(D) = 2\pi.
\]
However, by \cite[Theorem 4 p.189]{aleksandrovIntrinsicGeometrySurfaces2002},
\[
\tau_l(\gamma) = \tau_l(\gamma_1) + \tau_l(\gamma_2) + \pi - \tilde\theta_p + \pi - \tilde\theta_q,
\]
where $\tilde\theta_p$ and $\tilde\theta_q$ are the sector angles at $p,q$ on the inner side of the triangle.
Plugged into Gauss-Bonnet, we obtain
\[
\tau_l(\gamma_1) + \tau_l(\gamma_2) + \omega(D) = \tilde\theta_p + \tilde\theta_q.
\]
We distinguish two cases: first if anyone of the terms in the last equation is non-zero them
the left-hand side is non-positive as the right-hand side is non-negative.
This is a contradiction.
Otherwise, if all terms are zero, then $\omega\equiv 0$ on $D$.
Thus, the metric is Euclidean in $D$ and it contains a small Euclidean ball $B$ such that $\mu(B)>0$.
So $\omega(D)\leq\kappa\mu(D) < 0$.
This is a contradiction.
\end{proof}

\begin{proof}[Proof of theorem \ref{thm:BICAImpliesCAT}]
Let $(S,d)$ be a $\BICA(\kappa)$ surface and let $\kappa'>\kappa$.
Let $o\in S$ and let $P$ be a geodesically convex polygon containing $o$ and having a diameter less
than the constant $\ell$ given by theorem \ref{thm:CmBADiscImpliesHinge}.
We consider a triangle $T=[pqr]$ in $P$.
By the preceding lemma only branching phenomena can appear.
If the geodesics $[pq]$ and $[pr]$ branch in $p$ then the hinge comparison holds trivially because
\[
\measuredangle [p^q_r] = 0 \leq \modang^{\kappa'}(p^q_r).
\]
Otherwise, if the geodesics branches at the vertices $q$ or $r$ we can use the procedure given
in the proof of \cite[Theorem 10 p.216]{aleksandrovIntrinsicGeometrySurfaces2002} to deform, arbitrarily small,
$T$ into a triangle homeomorphic to a disc.
Then the hinge comparison holds
\[
\measuredangle [p^q_r] \leq \modang^{\kappa'}(p^q_r).
\]
Thus $P$ is a $\CAT(\kappa')$ space and then by using that $\CAT(\kappa')$ for every $\kappa'>\kappa$ implies $\CAT(\kappa)$,
$P$ is $\CAT(\kappa)$.
As $o$ was arbitrary, $S$ is locally $\CAT(\kappa)$.
\end{proof}

We close this paper with a few open questions that we have not been able to answer.
Even if the no-cusp assumption is very frequent in the theory, surfaces with cusps are natural objects so this would be interesting to extend
lemma \ref{lm:NotTooCuspsy} to know if $\BICB(\kappa)$ with $\kappa$ negative
implies the absence of cusp.
Also for curvature bounded above, is it possible to extend theorem \ref{thm:BICAImpliesCAT} to non-negative curvature?

\bibliography{bibliography}

\end{document}